%&amstex          
\input amstex\documentstyle{amsppt}  
\pagewidth{12.5cm}\pageheight{19cm}\magnification\magstep1
\topmatter
\title From Weyl groups to semisimple groups\endtitle
\author G. Lusztig\endauthor
\address{Department of Mathematics, M.I.T., Cambridge, MA 02139}\endaddress
\thanks{Supported by NSF grant DMS-1855773 and by a Simons Fellowship}\endthanks
\abstract{In this paper we show, using ideas from the theory of
total positivity, how a number of properties of a semisimple group
over the complex numbers can be presented purely in terms of the
Weyl group. We also describe some new connections of the theory of
canonical bases with total positivity.}\endabstract
\endtopmatter   
\document

\define\si{\sim}

\define\bin{\binom}
\define\op{\oplus}
   
\define\part{\partial}

\define\iy{\infty}
\define\m{\mapsto}
\define\do{\dots}

\define\sm{\smallmatrix}
\define\esm{\endsmallmatrix}
\define\sub{\subset}    

\define\T{\times}
\define\ti{\tilde}
\define\nl{\newline}
\redefine\i{^{-1}}

\define\un{\underline}

\define\a{\alpha}

\redefine\o{\omega}
\define\p{\pi}
\define\ph{\phi}

\define\r{\rho}

\define\th{\theta}
\define\k{\kappa}
\redefine\l{\lambda}

\define\x{\xi}

\define\Ps{\Psi}

\define\hh{\bold h}
\define\ii{\bold i}

\define\CC{\bold C}

\define\HH{\bold H}

\define\NN{\bold N}

\define\QQ{\bold Q}
\define\RR{\bold R}

\define\ZZ{\bold Z}

\define\cb{\Cal B}
\define\cc{\bold c}

\define\ci{\Cal I}

\define\fg{\frak g}

\define\fG{\frak G}

\define\tx{\ti x}

\define\tB{\ti B}

\define\sha{\sharp}

\define\bul{\bullet}

\define\che{\check}
\define\cha{\che{\a}}

\subhead 1\endsubhead
In this (partly expository) paper we show (using ideas from the theory of total positivity)
that many concepts related to a semisimple group $G$ over $\CC$ of simply laced type
can be presented purely in terms of the Weyl group.
This paper contains also a few new results. For example, we give a new characterization of the involution
$\ph$ studied in \cite{L97} in connection with the totally positive flag manifold.
In no.12 we show that the canonical basis \cite{L90}
of a finite dimensional irreducible representation of $G$ can be indexed by a set
which appears in the theory of total positivity (and whose definition involves the $\ZZ$-version of $\ph$).

In A3 we show that the totally positive flag manifold has something
close to a base point (a closed subset of dimension equal to the rank
of $G$).

In \S5, \S6 we state some new (conjectural) connections of
the theory of canonical bases with total positivity. (These
are verified  in some special cases in A4, A5.)

\subhead 2\endsubhead
We first define the Weyl group following Coxeter. (For simplicity we restrict ourselves to the simply laced case.)
Consider a finite connected graph with set of vertices $I'$ and with edges denoted by $i--j$ such that there
exists a function $h:I'@>>>\ZZ_{>0}$ with the following properties:

(1) for any $i\in I'$ we have $h(i)=(1/2)\sum_{j\in I';i--j}h(j)$ (harmonicity),

(2) $h(i)=1$ for some $i\in I'$.
\nl
Let $I$ be the graph obtained from $I'$ by removing one $i\in I'$ such that $h(i)=1$.
Coxeter has shown that the resulting graphs are exactly those that appear in the classification of
(simply laced) simple Lie algebras.

Here is an example of the graph $I'$ with the harmonic function $h$:
$$ \sm1&2&3&4&5&6&4&2\\
       {}&{}&{}&{}&{}&{3}&{}&{} \esm$$

The edges are pairs of numbers written next to each other. The graph $I$ (said to be of type $E_8$)
with the restriction of $h$ is       
$$\sm2&3&4&5&6&4&2\\
       {}&{}&{}&{}&{3}&{}&{} \esm$$
 In the rest of this paper the graph $I$ is fixed.
Let $E$ be the $\QQ$-vector space with basis $\{\cha_i;i\in I\}$. For $i\in I$ we define an automorphism
$s_i:E@>>>E$ by $\cha_j\m\cha_j-a_{ij}\cha_i$ where $a_{ij}$ is $2$ if $i=j$, is $-1$ if $i--j$, and is $0$
if $i\ne j$ don't form an edge of $I$. Let $W$ be the group of automorphisms of $E$ generated by
$\{s_i;i\in I\}$. This is the Weyl group. It is finite. For $w\in W$ we can write $w$ as a product of $s_i$;
the minimum number of factors in such a product is denoted by $|w|$. For example $|1|=0$; at the other extreme
there is a unique $w\in W$ for which $|w|$ is maximum; we denote it by $w_0$ and we set $\nu=|w_0|$.
Let $i\m i^!$ be the involution of $I$ such that $w_0s_iw_0=s_{i^!}$ for all $i\in I$.

\subhead 3\endsubhead
It is known that to our graph (or to $W$) corresponds a simply connected semisimple algebraic group $G$
over $\CC$.
Now $G$ has two important (unipotent) subgroups, $U^+,U^-$.
(For example to the graph with $I=\{i,j\}$ and with $i--j$ corresponds the algebraic group $SL_3(\CC)$
and $U^+,U^-$ is the group of upper triangular or lower triangular matrices with $1$ on diagonal; in this
case, $W$ is the symmetric group in $3$ letters.) We
would like to show how to construct $G$ from $W$. We will first try to construct $U^+,U^-$ from $W$.
A similar method applies to the full $G$ but this case will be only sketched.

\subhead 4\endsubhead
Let $U_{\ge0}$ be the semigroup with generators $\{i^a;i\in I,a\in\RR_{>0}\}$ and relations
(similar to those of a Coxeter group):

$i^ai^b=i^{a+b}$ for $i\in I$, $a,b$ in $\RR_{>0}$;

$i^aj^bi^c=j^{bc/(a+c)}i^{a+c}j^{ab/(a+c)}$ if $a_{ij}=-1$, $a,b,c$ in $\RR_{>0}$;

$i^aj^b=j^bi^a$ if $a_{ij}=0$, $a,b$ in $\RR_{>0}$.                                          

There is a unique semigroup anti-automorphism $\Ps:U_{\ge0}@>>>U_{\ge0}$ such that $\Ps(i^a)=i^a$ for all
$i\in I,a\in\RR_{>0}$. We have $\Ps^2=1$.

Let $\ci$ be the set of all sequences $\ii=(i_1,\do,i_\nu)$ in $I$ such that $w_0=s_{i_1}\do s_{i_\nu}$.
For $\ii\in\ci$ we define $\k_\ii:\RR_{>0}^\nu@>>>U_{\ge0}$ by
$$\cc=(c_1,\do,c_\nu)\m\ii^\cc:=i_1^{c_1}i_2^{c_2}\do i_\nu^{c_\nu}$$.
One can show that this map is injective and its image is independent of the choice of $\ii$.
We denote this image by $U_{>0}$. It is closed under multiplication in $U_{\ge0}$ and is stable under $\Ps$.

Let $\NN(Z_1,\do,Z_\nu)$ be the set of rational functions (coefficients in $\QQ$)
in the indeterminates $Z_1,\do,Z_\nu$ which are of the form $P(Z_1,\do,Z_\nu)/P'(Z_1,\do,Z_\nu)$ where
$P$ and $P'$ are (nonempty) sums of monomials in $Z_1,\do,Z_\nu$.

If $\ii\in\ci,\ii'\in \ci$, then by \cite{L94}, $\k_{\ii}\i\k_{\ii'}:\RR_{>0}^\nu@>>>\RR_{>0}^\nu$ is
of the form

(a) $(z_1,\do,z_\nu)\m(\p_1(z_1,\do,z_\nu),\do,\p_\nu(z_1,\do,z_\nu))$

where

(b) each $\p_1,\do,\p_\nu$ belongs to $\NN(Z_1,\do,Z_\nu)$.
\nl
It follows that $\k_{\ii}\i\k_{\ii'}$ can be regarded as a birational equivalence
$\CC^\nu-->\CC^\nu$. Let $O[\CC^\nu]$ be the algebra of regular functions $\CC^\nu@>>>\CC$ and let
$O(\CC^\nu)$ be the quotient field of this algebra. Now $\k_{\ii}\i\k_{\ii'}$ induces a field
isomorphism $O(\CC^\nu)@>>>O(\CC^\nu)$ denoted by $(\k_{\ii}\i\k_{\ii'})_*$ (it is given by sending an
element of $O(\CC^\nu)$ to its composition with $\k_{\ii}\i\k_{\ii'}$).
Let $\un O[U]$ be the set of all $(f_\ii)_{\ii\in\ci}$ where $f_\ii\in O(\CC^\nu)$ satisfy:

$(\k_\ii\i\k_{\ii'})_*(f_{\ii'})=f_\ii$ for any $\ii,\ii'$ in $\ci$;

$f_\ii\in O[\CC^\nu]$ for any $\ii\in\ci$.
\nl
This is a commutative  algebra in an obvious way.
The following result was conjectured in \cite{L19, \S6} and proved in \cite{FL21}.

(c) {\it $\un O[U]$ is the algebra of regular functions $O[U]$
on a unipotent algebraic group $U$ over $\CC$.}
\nl
Note that any element of $U_{>0}$ gives rise (via evaluation) to an algebra homomorphism $\un O[U]@>>>\CC$;
thus $U_{>0}$ can be regarded as a subset of $U$. The multiplication on $U$ extends that on $U_{>0}$ and
this defines it uniquely (by the requirement that it is regular). Now $U$ is the same as $U^+$ in no.3. 

\subhead 5\endsubhead
Let $O[U]_{\ge0}$ be the set of all $(f_\ii)_{\ii\in\ci}$ in $\un O[U]$ such that for any $\ii\in\ci$,
the function $\CC^\nu@>>>\CC$ given by $(c_1,\do,c_\nu)\m f_\ii(c_1,\do,c_\nu)$ is a polynomial
in $(c_1,\do,c_\nu)$ with coefficients in $\RR_{\ge0}$.
Note that $O[U]_{\ge0}$ is closed under addition, under multiplication and under scalar
multiplication by elements in $\RR_{\ge0}$ (but not under substraction).

One can also define $O[U]'_{\ge0}$ as the subset of $O[U]$ consisting of all $\RR_{\ge0}$-linear
combinations of the elements in the dual canonical basis \cite{L90} (at parameter $1$) of $O[U]$;
from the positivity properties of the canonical basis one can deduce that
$O[U]'_{\ge0}\sub O[U]_{\ge0}$. We conjecture that this inclusion
is an equality. (See A4 in the Appendix for a proof of this in a
special case.)

\subhead 6\endsubhead
The semigroup $\fG(\RR_{>0})$ defined in \cite{L19, 2.10} by generators
$i^a,\un i^a,-i^a$ ($i\in I,a\in\RR_{>0}$) and certain relations will be denoted here by $G_{\ge0}$.
We write $G_{>0}$ for the subset of $G_{\ge0}$ which in \cite{L19, 2.19} is denoted by
$\fG(\RR_{>0})_{w_0,-w_0}$; this is a sub-semigroup of $G_{\ge0}$.
Let $M=2\nu+|I|$. In \cite{L19, 2.13(b)} a family of bijections $\th_\hh:\RR_{>0}^M@>>>G_{>0}$ is described.
Here $\hh$ runs over a certain set of sequences with $M$ terms; we will take $\hh$ to be a sequence
of a special kind, that is either:

- the first $\nu$ terms form a sequence in $\ci$; the last $\nu$ terms form a sequence in $\ci$ (with the sign
$-$ attached) and the middle $|I|$ terms form a list of the elements of $I$ (underlined), or

-the first $\nu$ terms form a sequence in $\ci$ (with the sign $-$ attached); the last $\nu$ terms form a
sequence in $\ci$ and the middle $|I|$ terms form a list of the elements of $I$ (underlined).
\nl
These sequences form a finite set $\HH$.
The compositions $\th_\hh\i\th_{\hh'}:\RR_{>0}^M@>>>\RR_{>0}^M$ (with $\hh,\hh'$
in $\HH$) satisfy a property similar to 4(a),(b).
It follows that $\th_\hh\i\th_{\hh'}:\RR_{>0}^M@>>>\RR_{>0}^M$ can be regarded as a
birational equivalence $\CC^\nu\T(\CC^*)^{|I|}\T\CC^\nu-->\CC^\nu\T(\CC^*)^{|I|}\T\CC^\nu$.
Let $O[\CC^\nu\T(\CC^*)^{|I|}\T\CC^\nu]$ be the algebra of regular functions
$\CC^\nu\T(\CC^*)^{|I|}\T\CC^\nu@>>>\CC$ and let $O(\CC^\nu\T(\CC^*)^{|I|}\T\CC^\nu)$ be the quotient
field of this algebra. Now $\th_\hh\i\th_{\hh'}$ induces a field
isomorphism $O(\CC^\nu\T(\CC^*)^{|I|}\T\CC^\nu)@>>>O(\CC^\nu\T(\CC^*)^{|I|}\T\CC^\nu)$ denoted by
$(\th_\hh\i\th_{\hh'})_*$
(it is given by sending an element of $O(\CC^\nu\T(\CC^*)^{|I|}\T\CC^\nu)$ to its composition with
$\th_\hh\i\th_{\hh'}$). 
Let $\un O[G]$ be the set of all $(f_\hh)_{\hh\in\HH}$ where $f_\hh\in O(\CC^\nu\T(\CC^*)^{|I|}\T\CC^\nu)$
satisfy:

$(\th_\hh\i\th_{\hh'})_*(f_{\hh'})=f_\hh$ for any $\hh,\hh'$ in $\HH$;

$f_\hh\in O[\CC^\nu\T(\CC^*)^{|I|}\T\CC^\nu]$ for any $\hh\in\HH$.
\nl
This is a commutative  algebra in an obvious way.
The following result was conjectured in \cite{L19, \S6} and proved in \cite{FL21}.

(c) {\it $\un O[G]$ is the algebra $O[G]$ of regular functions on a semisimple simply connected
algebraic group $G$ over $\CC$.}
\nl
Note that any element of $G_{>0}$ gives rise (via evaluation) to an algebra homomorphism $\un O[G]@>>>\CC$;
thus $G_{>0}$ can be regarded as a subset of $G$. The multiplication on $G$ extends that on $G_{>0}$ and
this defines it uniquely (by the requirement that it is regular). Now this $G$ is the same as $G$ of no.3. 

Let $O[G]_{\ge0}$ be the set of all $(f_\hh)_{\hh\in\HH}$ in $\un O[G]$ such that for any $\hh\in\HH$,
the function $\CC^\nu\T(\CC^*)^{|I|}\T\CC^\nu@>>>\CC$ given by $(c_1,\do,c_M)\m f_\hh(c_1,\do,c_M)$ is an
$\RR_{\ge0}$-linear combination of functions $(c_1,\do,c_M)\m c_1^{k_1}\do c_M^{k_M}$
where $k_1,\do,k_M$ are integers of which the first $\nu$ and the last $\nu$ are $\ge0$.
Note that $O[G]_{\ge0}$ is closed under addition, under multiplication and under scalar
multiplication by elements in $\RR_{\ge0}$ (but not under substraction).

One can also define $O[G]'_{\ge0}\sub O[G]$ as the subset of $O[G]$ consisting of all $\RR_{\ge0}$-linear
combinations of the elements in the dual canonical basis (at parameter $1$) of $O[G]$. We conjecture that
$O[G]'_{\ge0}=O[G]_{\ge0}$. (See A5 in
the Appendix for a proof of this in a special case.)

\subhead 7\endsubhead
Now let $K$ be a semifield, that is a set with two operations: $+,\T$ such that $K$ is an abelian group
with respect to $\T$, an abelian semigroup with respect to $+$ and such that the distributivity law
$(a+b)c=ac+bc$ is satisfied. Here are three examples of semifields.

(i) $K=\RR_{>0}$ with the usual $+,\T$;

(ii) $K=\ZZ$ with the semifield structure in which the sum of $a,b$ is $\min(a,b)$ and the product of $a,b$
is $a+b$;

(iii) $K=\{1\}$ with $1+1=1,1\T1=1$.
\nl
If $\ii\in\ci,\ii'\in\ci$ and if we take $(z_1,\do,z_\nu)\in K^\nu$, then (in view of 4(b)), the
right hand side of 4(a) makes sense as an element of $K^\nu$, so that 4(a) defines a map
$(\k_\ii\i\k_{\ii'})_K:K^\nu@>>>K^\nu$ (which is inverse to $(\k_{\ii'}\i\k_\ii)_K$ hence is a bijection).

Let $U_K$ be the set of all $(\x_\ii)_{\ii\in\ci}$ where $\x_\ii\in K^\nu$ satisfy:

$(\k_\ii\i\k_{\ii'})_K(\x_{\ii'})=\x_\ii$ for any $\ii,\ii'$ in $\ci$.
\nl
Note that if $\ii\in\ci$, then $\ii$ defines a bijection $K^\nu@>\si>>U_K$ whose inverse is
$(\x_{\ii'})_{\ii'\in\ci}\m\x_\ii$. We denote this bijection by $\cc\m\ii^\cc$.
We have $U_{\RR_{>0}}=U_{>0}$. From the definitions we see that
the semigroup structure on $U_{>0}$ induces a semigroup structure on $U_K$.

There is a well defined involution $\Ps_K:U_K@>>>U_K$ given by
$(\x_\ii)_{\ii\in\ci}\m (\x'_\ii)_{\ii\in\ci}$
where for $\ii=(i_1,\do,i_\nu)\in\ci$ we have $\x'_\ii=(c_1,\do,c_\nu)$
whenever $\x_{i_\nu,\do,i_1}=(c_\nu,\do,c_1)$. This is an anti-automorphism of $U_K$;
when $K=\RR_{>0}$, it coincides with the restriction of $\Ps:U_{\ge0}@>>>U_{\ge0}$ (in no.4)
to $U_{>0}$.

\subhead 8\endsubhead
For $p=(p_i)\in K^I$ we define a bijection $S_p:U_K@>>>U_K$ by
$$(\x_\ii)_{\ii\in\ci}\m(\x'_\ii)_{\ii\in\ci}$$
where $\x'_\ii$ is obtained by multiplying $\x_\ii$ by $(p_{i_1},p_{i_2},\do,p_{i_\nu})$
(component by component in $K$) where $\ii=(i_1,\do,i_\nu)$. See \cite{L19, 4.3(a)}.
From the definitions we have $S_p\Ps_K=\Ps_KS_p$.

We define an element $u(1)\in U_K$ as follows. Let $\ii\in\ci$. For $k\in[1,\nu]$ we have
$$s_{i_1}s_{i_2}\do s_{i_{k-1}}(\cha_{i_k})=\sum_{i\in I}r_{i,k}\cha_i$$
where $r_{i,k}\in\NN$. Let $r'_k=\sum_{i\in I}r_{i,k}\in\ZZ_{>0}$. Let
$\cc=(r'_1,r'_2,\do,r'_\nu)\in K^\nu$. As in \cite{L94, 11.2} we see that $\ii^\cc\in U_K$ is
independent of the choice of $\ii$; we denote it by $u(1)$.
We define $q=(q_i)\in K^I$ by $q_i=\sum_{k\in[1,\nu]}r_{i,k}\in\ZZ_{>0}\sub K$; this is also
independent of the choice of $\ii$.

We define an imbedding $K^I@>>>U_K$ by $p\m u(p):=S_p(u(1))$.

\subhead 9\endsubhead
For $\cc=(c_1,c_2,\do,c_\nu)\in K^\nu$ and $c\in K$ we set
${}_c\cc=(cc_1,c_2,\do,c_\nu)\in K^\nu$.
For any $i\in I,c\in K$ there is a unique bijection $T_{i,c}:U_K@>>>U_K$ such that
$T_{i,c}(\ii^\cc)=\ii^{{}_c\cc}$ for some/any $\ii=(i_1,\do,i_\nu)\in\ci$ such that $i_1=i$ and any
$\cc\in K^\nu$ (see \cite{L97, 2.3} for the case $K=\ZZ$ and \cite{L19, 2.16} for a general $K$).

We regard $U_K$ as the set of vertices of a graph in which $u\ne u'$ are joined if $u'=T_{i,c}(u)$
for some $i\in I,c\in K$. We have the following result.

(a) {\it If $K=\RR_{>0}$, this graph is connected.}
\nl
The proof is given in the Appendix, see A2. (An analogous result in which $\RR_{>0}$ is replaced by
$\ZZ$ in 7(ii) appears in \cite{L97, 2.8}.)

\subhead 10\endsubhead
We have the following result. (An analogous result in which $\RR_{>0}$ is replaced by $\ZZ$ in 7(ii) appears
in \cite{L97, 2.9}.)

(a) {\it There is a  unique bijection $\ph:U_{>0}@>>>U_{>0}$ such that}

(i) $T_{i,c}\ph=\ph T_{i^!,c\i}$ for all $i\in I,c\in\RR_{>0}$,

(ii) $\ph(u(1))=u(q\i)$.
\nl
The existence is proved in the Appendix, see A1. The uniqueness of $\ph$ follows from 9(a). 

For example, if $I=\{i\}$, then $\ph$ is given by $i^c\m i^{c\i}$; if $I=\{i,j\}$ with $i--j$, then
$\ph$ is given by 
$$i^aj^bi^c\m i^{a/c(a+c)}j^{(a+c)/ab}i^{1/(a+c)}=j^{c/ab}i^{1/c}j^{1/b}.$$

\subhead 11\endsubhead
The results in this subsection are based on the identification of $\ph$ in no.10 with the bijection
with the same name in \cite{L97} (see the proof in A1).
For $p\in\RR_{>0}^I$ we have $S_p\ph=\ph S_{p\i}$. (See \cite{L19, 4.3(d)}.) It follows that
$\ph(u(p))=u(q\i p\i)$. We have $\ph^2=1$ (see \cite{L19, 4.1}).

From \cite{L97, 3.4} it follows that if $\ii\in\ci,\ii'\in\ci$, then
$\k_\ii\i\ph\k_{\ii'}:\RR_{>0}^\nu@>>>\RR_{>0}^\nu$ is of the form
$$(z_1,\do,z_\nu)\m(\r_1(z_1,\do,z_\nu)),\do,\r_\nu(z_1,\do,z_\nu))$$
where each $\r_1,\do,\r_\nu$ belongs to $\NN(Z_1,\do,Z_\nu)$.
It follows that if $K$ is a semifield, this map gives rise to a map $K^\nu@>>>K^\nu$ which can be viewed as a 
bijection $U_K@>>>U_K$ (denoted by $\ph_K$) which does not depend on the choice of $\ii,\ii'$.

We have $\ph_K^2=1$. Since for $p\in K^I$ we have $S_p\ph_K=\ph_K S_{p\i}$, we see that $(S_p\ph_K)^2=1$.

We set $\ph'_K=\Ps_K\ph_K\Ps_K:U_K@>>>U_K$, ($\Ps_K$ as in no.7). 
For $p\in K^I$ we have $S_p\ph'_K=\ph'_KS_{p\i}$, hence $(S_p\ph'_K)^2=1$.

\subhead 12\endsubhead
We now assume that $K=\ZZ$ is as in 7(ii).
For $i\in I$ we define $z_i:U_\ZZ@>>>\ZZ$ by $z_i((\x_\ii)_{\ii\in\ci})=c_\nu$
where $\x_\ii=(c_1,\do,c_{\nu-1},c_\nu)$ is defined in terms of $\ii=(i_1,\do,i_\nu)\in\ci$ such that
$i_\nu=i$. (It is easy to see that $z_i$ is well defined).

Let $U_\NN$ be the set of all $(\x_\ii)_{\ii\in\ci}\in U_\ZZ$ such that for any $\ii\in\ci$,
$\x_\ii$ is not only in $\ZZ^\nu$ but actually in $\NN^\nu$. (Note that, if $\x_\ii\in\NN^\nu$ for some $\ii$,
then $\x_\ii\in\NN^\nu$ for all $\ii$, see \cite{L19, 2.14}.)

Let $\l=(\l_i)_{i\in I}\in\NN^I$. Now $S_\l:U_\ZZ@>>>U_\ZZ$ is defined as in no.8 since $\l\in\ZZ^I$. It is
known that $\l$ indexes a finite dimensional irreducible representation $V_\l$ of $G$ (in no.3) with a
canonical basis \cite{L90} in natural bijection with
$$U_{\NN,\l}:=\{x\in U_\NN;z_i(x)\le\l_i \text{ for all }i\in I\}.$$
Let
$$U'_{\NN,\l}=\{x\in U_\NN;S_\l\ph_\ZZ(x)\in U_\NN\}.$$
We shall prove the following result which appears as a conjecture in \cite{L19, 8.2(b)}.

(a) {\it We have $U_{\NN,\l}=U'_{\NN,\l}$.}
\nl
We shall write $\ph$ instead of $\ph_\ZZ$.
From \cite{L97, 4.9} (see also the errata in \cite{L17, \#130}) we have

(b) $S_\l\ph:U_{\NN,\l}@>\si>>U_{\NN,\l}$.
\nl
Since $U_{\NN,\l}\sub U_\NN$, it follows that $U_{\NN,\l}\sub U'_{\NN,\l}$. Now let $x\in U'_{\NN,\l}$. 
We set $\tx=S_\l\ph(x)\in U_\NN$. Define $\mu=(\mu_i)\in\NN^I$ by $\mu_i=z_i(\tx)$.  We have
$\tx\in U_{\NN,\mu}$. By (b) (for $\mu$ instead of $\l$) we have $\tx=S_\mu\ph(y)$ for some
$y\in U_{\NN,\mu}$. Now $S_\l\ph(x)=S_\mu\ph(y)$ hence $\ph(S_{-\l}(x))=\ph(S_{-\mu}(y))$ so that
$S_{-\l}(x)=S_{-\mu}(y)$ and $x=S_{\l-\mu}(y)$. We have $z_i(x)=\l_i-\mu_i+z_i(y)$.
Since $y\in U_{\NN,\mu}$, we have $z_i(y)\le\mu_i$, hence $\l_i-\mu_i+z_i(y)\le\l_i$ and
$z_i(x)\le\l_i$. We see that $x\in U_{\NN,\l}$. This proves (a).

Let 
$$U''_{\NN,\l}=\{x\in U_\NN;S_\l\ph'_\ZZ(x)\in U_\NN\}.$$
We have a bijection $U'_{\NN,\l}@>>>U''_{\NN,\l}$ given by $x\m\Ps_\ZZ(x)$. (Note that
$\Ps_\ZZ$ preserves $U_\NN$.) Thus $U'_{\NN,\l}$ can be identified via $\Ps_\ZZ$ with
$U''_{\NN,\l}$.

In particular we have

(c) $\dim(V_\l)=\sha(U'_{\NN,\l})=\sha(U''_{\NN,\l})$.
\nl
Earlier formulas for $\dim(V_\l)$ were given by Weyl (as a quotient of two positive integers) and by  Kostant
(as a difference of two integers); in both of these formulas the result was not
obviously a positive integer. On the other hand, last expression in (c) is either a positive
integer or $\iy$. Combining any two of these three formulas shows that the result is a positive integer.

\head Appendix\endhead
\subhead A1\endsubhead
In this subsection we give a proof of the existence part of 10(a).
Let $G$ (over $\CC$) be as in no.3. Let $\fg$ be the Lie algebra of $G$.
We assume given a maximal torus $T$ of $G$ and a pair $B^+,B^-$ of opposed Borel subgroups
of $G$ containing $T$, with unipotent radicals $U^+,U^-$. For $i\in I$ we consider homomorphisms
$x_i:\CC@>>>U^+,y_i:\CC@>>>U^-$ such that $(T,B^+,B^-,x_i,y_i;i\in I)$ is a pinning for $G$.
Define $e_i,f_i$ in $\fg$ by $\exp(ae_i)=x_i(a),\exp(af_i)=y_i(a)$ for all $a\in\CC$. Let $h_i=[e_i,f_i]$.
There is a unique semigroup imbedding $U_{\ge0}@>>>U^+$, $u\m u^+$, given by $i^c\m x_i(c)$ for any
$i\in I$ and any $c\in\RR_{>0}$. There is a unique semigroup imbedding $U_{\ge0}@>>>U^-$, $u\m u^-$, given by
$i^c\m y_i(c)$ for any $i\in I$ and any $c\in\RR_{>0}$.

By \cite{L97, 3.3} there exists a unique bijection
$\ph:U_{>0}@>>>U_{>0}$ such that
$$(\ph(u)^-)\i B^+\ph(u)^-=(u^+)\i B^-u^+$$
for all $u\in U_{>0}$.
As stated in \cite{L19, 4.2(a)}, this $\ph$ satisfies 10(a)(i). (This follows immediately from
\cite{L97, Lemma 3.6}). To verify that this $\ph$ satisfies 10(a)(ii) it is enough to show that
$$(\ph(u(1))^-)\i B^+\ph(u(1))^-=(u(q\i)^-)\i B^+u(q\i)^-$$
or equivalently that
$$(u(1)^+)\i B^-u(1)^+=(u(q\i)^-)\i B^+u(q\i)^-.$$
According to the conjecture \cite{L94, 11.4(a)} (with all $p_i=1$), proved in \cite{FL97}, we have
$u(1)^+=\exp(\sum_{i\in I}q_ie_i)$ (with $q_i$ as in no.8). The same proof applied 
with $e_i,f_i$ replaced by $q_i\i e_i,q_if_i$ shows that $u(q\i)^+=\exp(\sum_{i\in I}e_i)$ hence
$u(q\i)^-=\exp(\sum_{i\in I}f_i)$. Thus we are reduced to proving
$${}^{\exp(-\sum_{i\in I}q_ie_i)}B^-={}^{\exp(-\sum_{i\in I}f_i)}B^+$$
or, setting $\o=\sum_{i\in I}q_ie_i,\o'=\sum_{i\in I}f_i$, that
$$\exp(\o')\exp(-\o)B^-\exp(\o)\exp(-\o')=B^+.\tag a$$
Here for $g\in G,B\in\cb$ we write ${}^gB$ instead of $gBg\i$.
We have $[\o,\o']=\sum_{i\in I}q_ih_i$. We show that $[[\o,\o'],\o]=2\o$, $[[\o,\o'],\o']=-2\o'$ or
equivalently that
$$\sum_{i,j}q_iq_ja_{ij}e_j=2\sum_{i\in I}q_ie_i,-\sum_{i,j}q_ia_{ij}f_j=-2\sum_{i\in I}f_i$$
or equivalently that $\sum_iq_ia_{ij}=2$. This follows  from the definition of $q_i$.
We see that $\o,\o',[\o,\o']$ is an $sl_2$-triple. Hence the subgroup of $G$ with Lie algebra
$\CC\o\op\CC\o'\op\CC[\o,\o']$ is isomorphic to $SL_2(\CC)$ or $PGL_2(\CC)$. By a property of $SL_2(\CC)$, we
see that $Ad(\exp(\o')\exp(-\o)\exp(\o'))$ maps the line spanned by $\o'$ to the line spanned by $\o$.
Since $\o$ is regular nilpotent in $Lie(B^+)$, we see that $B^+$ is the unique Borel subgroup whose Lie
algebra contains $\o$; similarly $B^-$ is the unique Borel subgroup whose Lie algebra contains $\o'$. It
follows that conjugation by $\exp(\o')\exp(-\o)\exp(\o')$ takes $B^-$ to $B^+$. Since $\exp(\o')\in B^-$, it
follows that (a) holds. This completes the proof of the existence part of 10(a).

\subhead A2\endsubhead
In this subsection we give a proof of 9(a). We use notation in A1. Let $\cb$ be the variety of
Borel subgroups of $G$. 
If $B\in\cb,B'\in\cb$ we denote by $pos(B,B')$ the relative position of $B,B'$
(an element of $W$, the Weyl group of $G$).
Let $G(\RR)$ be the group of real points of $G$ defined by the pinning and let $\cb(\RR)$ be the orbit of
$B^+$ (or $B^-$) under the adjoint action of $G(\RR)$.

If $i\in I$, an $i$-circle in $\cb(\RR)$ is a subset of
$\cb(\RR)$ of the form
$\{B\}\cup\{B'\in\cb(\RR);pos(B,B')=s_i\}$ for some
$B\in\cb(\RR)$.

Let
$$\cb_{>0}=\{{}^{u^+}B^-;u\in U_{>0}\}=\{{}^{u^-}B^+;u\in U_{>0}\}.$$
(The last equality follows from \cite{L94}.) This is an open subset of $\cb(\RR)$.

Let $u\ne u'$ in $U_{>0}$ and $i\in I$ be such that

(a) $pos({}^{u^-}B^+,{}^{u'{}^-}B^+)=s_i$.
\nl
We show that

(b) $\Ps(u')=T_{i,a}\Ps(u)$ for some $a\in\RR_{>0}$, ($\Ps$ as in no.4).              
\nl
Now $(u\i u')^-$ is in the intersection of $U^-$ with the parabolic
subgroup generated by $B^+$ and by $y_i(\CC)$ hence is in $y_i(\CC)$. It follows that $(u')^-=u^-y_i(c)$
for some $c\in\CC$ (which is necessarily in $\RR-\{0\}$). If $c>0$ we deduce $\Ps(u')=i^c\Ps(u)$, from which
(b) follows immediately. If $c<0$ then $u^-=(u')^-y_i(-c)$ with $-c>0$ which implies 
$\Ps(u)=i^{-c}\Ps(u')$ and again (b) follows.
(Conversely, it is easy to show that if (b) holds then (a) holds.)
Note also that if $c>0$ then $\{{}^{u^-y_i(c')}B^+;0\le c'\le c\}$ is contained in $\cb_{>0}$ and in an
$i$-circle and it contains both ${}^{u^-}B^+,{}^{u'{}^-}B^+$; if $c<0$ then
$\{{}^{u'{}^-y_i(c')}B^+;0\le c'\le -c\}$ is contained in $\cb_{>0}$ and in an
$i$-circle and it contains both ${}^{u^-}B^+,{}^{u'{}^-}B^+$. We see that the intersection of any $i$-circle
in $\cb(\RR)$ with $\cb_{>0}$ is either empty or connected. (If it is nonempty, this intersection is called
a half $i$-circle.)

We define a graph structure on $\cb_{>0}$ in which $B\ne B'$ in $\cb_{>0}$ are joined if for some $i\in I$,
we have $pos(B,B')=s_i$ or equivalently if $B,B'$ belong to the same half $i$-circle.
We have the following result.

(c) {\it This graph on $\cb_{>0}$ is connected.}
\nl
Since (a) implies (b), we see that to prove 9(a) it is enough to prove (c).

We first verify the following statement.

(d) Let $\ii=(i_1,\do,i_\nu)\in\ci$. If $B\in\cb_{>0}$, then for any $B'\in\cb_{>0}$ sufficiently close to
$B^+$ there is a unique sequence $B'=B_0,B_1,\do,B_\nu=B$ in $\cb_{>0}$ such that
$pos(B_0,B_1)=s_{i_1},pos(B_1,B_2)=s_{i_2},\do,pos(B_{\nu-1},B_\nu)=s_{i_\nu}$.
\nl
We have $B={}^{u^-}B^+$ where $u\in U_{>0}$.

If  $u'\in U_{>0}$ is such that $u'{}^-$ is sufficiently close to $1$ then:

(i) $pos(B,{}^{(u')^-}B^+)=w_0$ and

(ii) $(u'{}\i u)^-\in U^-_{>0}$.
\nl
Indeed, for $u'{}^-$ close to $1$, (ii) holds since $u^-\in U^-_{>0}$ and
$U^-_{>0}$ is open in the group of real points of $U^-$. Also, we have $pos(B,B^+)=w_0$ so that
for $u'{}^-$ close to $1$, (i) holds (we use that ${}^{(u')^-}B^+$ is close to $B^+$ hence
is contained in the open subset $\{B_1\in\cb;pos(B,B_1)=w_0\}$ of $\cb$).

We write $(u'{}\i u)^-=y_{i_1}(c_1)\do y_{i_\nu}(c_\nu)$ where $c_1,\do,c_\nu$ are in $\RR_{>0}$. We set

$B_0={}^{(u')^-}B^+$, $B_1={}^{(u')^-y_{i_1}(c_1)}B^+$, $B_2={}^{(u')^-y_{i_1}(c_1)y_{i_2}(c_2)}B^+$, $\do$,

$B_\nu={}^{(u')^-y_{i_1}(c_1)\do y_{i_\nu}(c_\nu)}B^+$.   

For $k=0,1,2,\do,\nu$ we have $u'y_{i_1}(c_1)\do y_{i_k}(c_k)\in U^-_{>0}U^-_{\ge0}\sub U^-_{>0}$
hence $B_0,B_1,\do,B_\nu$ are in $\cb_{>0}$. Note that $B_\nu={}^{u^-}B^+=B$. We have

$pos(B_0,B_1)=pos(B^+,{}^{y_{i_1}(c_1)}B^+)=s_1$, $pos(B_1,B_2)=pos(B^+,{}^{y_{i_2}(c_2)}B^+)=s_2$, $\do$,

$pos(B_{\nu-1},B_\nu)=pos(B^+,{}^{y_{i_\nu}(c_\nu)}B^+)=s_\nu$.
\nl
This proves the existence in (d). The uniqueness is obvious.

We now prove (c). Let $B\ne\tB$ in $\cb_{>0}$. If $B'\in\cb_{>0}$ is sufficiently close to $B^+$ then
by (d), $B$ can be joined with $B'$ through a sequence of edges of our graph and $\tB$
can be joined with $B'$ through a sequence of edges of our graph. Thus $B,\tB$ are in the same connected component of
our graph. This completes the proof of (c) hence that of 9(a).

\subhead A3\endsubhead
Let $p\in\RR_{>0}^I$. From the equality $\ph(u(p))=u(q\i p\i)$ (see no.11) and
the equality $(\ph(u)^-)\i B^+\ph(u)^-=(u^+)\i B^-u^+$ for $u\in U_{>0}$ we deduce

$(u(p)^+)\i B^-u(p)^+=(u(q\i p\i)^-)\i B^+u(q\i p\i)^-$.
\nl
Applying the antiautomorphism of $G$ which keeps each $x_i(a),y_i(a)$ fixed
(hence keeps $u(p)^+,u(p)^-$ fixed), we deduce

$u(p)^+B^-(u(p)^+)\i=u(q\i p\i)^-B^+(u(q\i p\i)^-)\i$.
\nl
It follows that
$$\{{}^{u(p)^-}B^+;p\in \RR_{>0}^I\}=\{{}^{u(p)^+}B^-;p\in \RR_{>0}^I\}.$$
The two sides of this equality form a (closed) subset $\cb_{>0}^\bul$ of $\cb_{>0}$
which is a single orbit under a (free) $\RR_{>0}^I$-action on $\cb_{>0}$.
This subset is the closest we can come to having a base point of $\cb_{>0}$.

\subhead A4\endsubhead
In this subsection we assume that $I=\{i,j\}$, $i--j$. In this case the vector space $\un O[U]$ in no.4
consists of all pairs $[\Pi;\Pi']$ where $\Pi$ and $\Pi'$ are polynomials in the indeterminates $a,b,c$
with coefficients in $\CC$ and we have
$$\Pi(bc/(a+c),a+c,ab/(a+c))=\Pi'(a,b,c).$$
The following are examples of pairs in $\un O[U]$:
$$\align&f_{i,j,k}=[a^ib^jc^k;((bc/(a+c))^i(a+c)^j(ab/(a+c))^k]\\&
=[a^ib^jc^k;\sum_{h',h''\text{ in }\NN;h'+h''=j-i-k}\bin{j-i-k}{h'}a^{h'+k}b^{i+k}c^{h''+i}]\tag a\endalign$$
with $i,j,k$ in $\NN$, $j\ge i+k$,
$$\align&f'_{i,j,k}=[((bc/(a+c))^i(a+c)^j(ab/(a+c))^k;a^ib^jc^k]\\&
=[\sum_{h',h''\text{ in }\NN;h'+h''=j-i-k}\bin{j-i-k}{h'}a^{h'+k}b^{i+k}c^{h''+i};a^ib^jc^k]\tag b\endalign$$
with $i,j,k$ in $\NN$, $j\ge i+k$.

(c) These pairs are distinct except for the equality $f_{i,j,k}=f'_{k,j,i}$ when $j=i+k$.
\nl
We show that they are linearly independent. Indeed, $\un O[U]=O[U]$ is a direct sum of weight spaces
coming from a $(\CC^*)^2$-action on $U$ and indexed by pairs $(m,n)\in\NN^2$.
The pair $f_{i,j,k}$ is in the weight space indexed by $(i+k,j)$ and the pair
$f'_{i,j,k}$ is in the weight space indexed by $(j,i+k)$. Hence the pairs (a),(b) in a given weight
space are all of type (a) or all of type (b). But the pairs of type (a) are linearly independent
(since their first component are clearly linearly independent) and the pairs of type (b) are linearly
independent (since their second component are clearly linearly independent).
Thus the pairs (a),(b) (with the identification (c)) are linearly independent.
We can verify that the number of pairs (a),(b) in a given weight space is equal to the known dimension
of that weight space. It follows that the pairs (a),(b) (with the identification (c))
form a basis of $\un O[U]=O[U]$.

Let $O[U]''_{\ge0}$ be the set of $\RR_{\ge0}$-linear combinations of pairs in (a),(b).
Since each pair in (a),(b) belongs to $O[U]_{\ge0}$, we have $O[U]''_{\ge0}\sub O[U]_{\ge0}$.  
Conversely consider a pair $[x;x']\in O[U]_{\ge0}$. Then the projections of $[x;x']$ to the
various weight spaces of $\un O[U]$ are also in $O[U]_{\ge0}$. Hence to prove that
$[x;x']\in O[U]''_{\ge0}$ we can assume that $[x;x']$ is in a weight space of $\un O[U]$
(and in $O[U]_{\ge0}$). Thus we have either

(i) $[x;x']=\sum_{j\ge i+k}c_{i,j,k}f_{i,j,k}$ or

(ii) $[x;x']=\sum_{j\ge i+k}c'_{i,j,k}f'_{i,j,k}$
\nl
with $c_{i,j,k}\in\CC$, $c'_{i,j,k}\in\CC$. If (i) holds then
$$x=\sum_{j\ge i+k}c_{i,j,k}a^ib^jc^k$$ and since
$[x;x']\in O[U]_{\ge0}$ we must have $c_{i,j,k}\in\RR_{\ge0}$.
If (ii) holds then $$x'=\sum_{j\ge i+k}c'_{i,j,k}a^ib^jc^k$$ and since
$[x;x']\in O[U]_{\ge0}$ we must have $c'_{i,j,k}\in\RR_{\ge0}$.
We see than in either case we have $[x;x']\in O[U]''_{\ge0}$.
Thus we have $O[U]_{\ge0}=O[U]''_{\ge0}$.

One can show that under the identification $\un O[U]=O[U]$ in 4(c) the pairs (a),(b)
form precisely the dual canonical basis \cite{L90} (with parameter $1$).
It follows that in this case we have $O[U]_{\ge0}=O[U]'_{\ge0}$ (see no.5).
We see that (at least in this case) the dual canonical basis of $O[U]$ can be recovered (up to
multiplication by scalars in $\RR_{>0}$) without using the theory of quantum groups and without
intersection cohomology: it is the only basis of $O[U]$ (up to multiplication by scalars in $\RR_{>0}$) such
that the set of $\RR_{\ge0}$-linear combination of its elements is exactly $O[U]_{\ge0}$.

\subhead A5\endsubhead
In this subsection we assume that $I=\{i\}$. In this case the vector space $\un O[G]$ in no.6
consists of all pairs $[\Pi;\Pi']$ where $\Pi$ and $\Pi'$ are polynomials in the
indeterminates $a,b,b\i,c$ with coefficients in $\CC$ and we have
$$\Pi(a,b\i,c)=\Pi'(c/(ac+b^2),b/(ac+b^2),a/(ac+b^2)).$$
The following are examples of pairs in $\un O[G]$:
$$\align&g_{i,j,k}=[a^ib^{-j}c^k;(c/(ac+b^2))^i(b/(ac+b^2))^j(a/(ac+b^2))^k]\\&
=[a^ib^{-j}c^k;
\sum_{h',h''\text{ in }\NN;h'+h''=-i-j-k}\bin{-i-j-k}{h'}a^{k+h'}b^{j+2h''}c^{i+h'}]\tag a\endalign$$
with $i,k$ in $\NN$, $j\in\ZZ$, $-i-j-k\ge0$,
$$\align&g'_{i,j,k}=[(c/(ac+b^{-2}))^i(b\i/(ac+b^{-2}))^j(a/(ac+b^{-2}))^k;a^ib^jc^k]\\&
=[\sum_{h',h''\text{ in }\NN;h'+h''=-i-j-k}a^kc^ib^{2i+2k+2j-j}(ab^2c+1)^{-i-j-k};a^ib^jc^k]\tag b\endalign$$

with $i,k$ in $\NN$, $j\in\ZZ$, $-i-j-k\ge0$.

(c) These pairs are distinct except for the equality $g_{i,j,k}=g'_{k,j,i}$ when $i+j+k=0$.
\nl
One can verify that the pairs (a),(b) (with the identification (c)) form precisely the dual canonical
basis (with parameter $1$) of $\un O[G]=O[G]$; these pairs are clearly in $O[G]_{\ge0}$
and one can verify that $O[G]_{\ge0}=O[G]'_{\ge0}$.

\enddocument